\documentclass[12pt]{article}
\usepackage[utf8]{inputenc}
\usepackage{amsmath}
\usepackage{amsthm}
\usepackage{amscd}
\usepackage{comment}
\usepackage[english]{babel}
\usepackage{tikz}
\usepackage{color,graphicx,amssymb,tikz-cd,wrapfig}
\usepackage{amsmath,amsthm,amsfonts,amssymb,stmaryrd}
\usepackage{bm}
\usepackage[rightcaption]{sidecap}
\usepackage{geometry}
\usepackage{hyperref}
\usepackage[figurename=Fig.]{caption}

\usepackage[T1]{fontenc}
\usepackage{combelow}
\usepackage{newunicodechar}

\usepackage{lipsum}   

\usepackage{xcolor}

\usepackage{enumerate}
\usepackage{pb-diagram}
\usepackage{graphicx} 
\usepackage[all]{xy}
\usepackage{appendix}
\usetikzlibrary{intersections,calc,arrows.meta}

\title{Relative norm averaging property and selflessness of almost periodic actions}

\author{Ryo Toyota\thanks{Department of Mathematical Sciences, University of Copenhagen, Universitetsparken 5, DK-2100, Copenhagen \O, Denmark}}%
\date{ }

\newtheorem{thm}{Theorem}[section]
\newtheorem{dfn}[thm]{Definition}
\newtheorem{lem}[thm]{Lemma}
\newtheorem{cor}[thm]{Corollary}
\newtheorem{rmk}[thm]{Remark}

\newtheorem{ex}[thm]{Example}

\newcommand{\CC}{ {\mathbb C} }
\newcommand{\RR}{ {\mathbb R} }
\newcommand{\ZZ}{ {\mathbb Z} }

\newcommand{\FF}{{\mathbb F}}

\newcommand{\id}{\text{id}}

\newcommand{\aut}{\text{Aut}}
\newcommand{\conv}{\text{conv}}
\newcommand{\Tr}{\text{Tr}}
\newcommand{\GL}{\text{GL}}

\begin{document}

\maketitle

\begin{abstract}
Let $D$ be a unital separable $C^*$-algebra with a faithful state $\sigma$ and $G$ be a discrete group.
We introduce a property called relative norm averaging property (RNA) for almost periodic $\sigma$-preserving actions by $G$ on $(D,\sigma)$.

We prove that if an almost periodic $\sigma$-preserving action $\alpha$ satisfies that $\overline{\alpha(G)} \subseteq \aut(D,\sigma)$ is connected and nontrivial, then $\alpha$ satisfies (RNA). As an application, we prove that if $G$ has property PHP of Ozawa and $\alpha:G\curvearrowright (D,\sigma)$ is an ergodic (i.e. $D^G=\mathbb{C}$) almost periodic $\sigma$-preserving action with $\overline{\alpha(G)}$ being connected, then $D\rtimes_r G$ is selfless. This includes certain ergodic actions by groups with property (T) and paradoxical isometric actions of $\mathbb{F}_2$ on the sphere.

In non-ergodic case, we prove that the relative norm averaging property implies that $D\rtimes_r G$ satisfies the $D^G$-valued selfless, if $G$ satisfies PHP.
\end{abstract}

\section{Introduction}
Stable rank, introduced in \cite{RieffelStable}, is a noncommutative notion of dimension. A unital $C^*$-algebra $A$ is said to have stable rank one if its invertible elements are dense in $A$. This approximation property plays a fundamental role in the structure theory of $C^*$-algebras and substantial progress has been made in establishing stable rank one in a variety of settings \cite{DykemaHaagerupRordam,RordamStableRealRank,GerasimovaOsin}.

Determining when reduced crossed products have stable rank one is a natural problem connecting dynamics with the structure of $C^*$-algebras. For free minimal $\ZZ^d$-actions, stable rank one was established by Phillips for Cantor systems \cite{PhillipsCrossed} and by Li and Niu for arbitrary compact Hausdorff spaces \cite{LiNiuStableRank}. Suzuki also proved stable rank one for the crossed products associated with universal minimal systems of countable amenable groups \cite{SuzukiAlmostFinite}. 

In the nonnuclear setting, recently Bell, Geffen, and Kerr introduced the concept of squarely divisibility and proved that for a topologically free, minimal, square divisible action $G\curvearrowright X$, which admits an invariant Borel probability measure, the stable rank of $C(X)\rtimes_r G$ is one \cite{BellGeffenKerrStable}. Also they prove that stable rank one is generic among minimal, topologically free, topologically weakly mixing Cantor actions of finitely generated free groups that admit an invariant Borel probability measure.

In this paper, we study almost periodic actions (Definition~\ref{almost periodic}), motivated by isometric actions on compact Riemannian manifolds,  including certain
actions admitting paradoxical decompositions and actions with
spectral gap (see Example~\ref{example}). Our approach uses selflessness, introduced by Robert \cite{RobertSelfless}, which in the tracial setting implies stable rank one. 

We establish a relative norm averaging property for almost periodic actions whose compact closure is connected. For ergodic actions of groups with Ozawa’s property PHP, this yields complete selflessness of the reduced crossed product. We also prove an operator-valued selflessness theorem for nonergodic actions.

\vspace{.2in}

In this paper, let us denote by $D$ a unital separable $C^*$-algebra with a faithful state $\sigma$, by $G$ a countable discrete group with a $\sigma$-preserving action $\alpha:G \to \aut(D,\sigma)$. We also denote the canonical conditional expectation by $E:D\rtimes_r G \to D$.
\begin{dfn}\label{almost periodic}
    We say that the action $\alpha:G \to \aut(D)$ is almost periodic if $K:=\overline{\alpha(G)}\subseteq \aut(D)$ is compact with respect to the point-norm topology. In this case, we denote by $P_K$ the average with respect to the Haar measure $m_K$ on $K$, i.e. 
    \begin{align*}
        P_K:D\to D^G;\; a\mapsto \int_{K}\kappa(a) dm_K(\kappa).
    \end{align*}
\end{dfn}

Note that if $\overline{\alpha(G)}$ is compact, a pointwise limit $\alpha:D \to D$ of $\{\alpha_{g_i}\}_{g_i}$ is again in $\aut(D)$. We must check the surjectivity of $\alpha$. By the compactness of $K$,  there exists a limit point $\beta$ of $\{\alpha_{g_i}^{-1}\}$. Then we have $\alpha\circ \beta=\id$ and $\alpha$ must be surjective. 

The central property of this paper is a relative version of Bell’s norm averaging property \cite{BellSelflessNonnuclear}.

\begin{dfn}
    We say that an almost periodic action $\alpha:G \to \aut(D)$ with $K:=\overline{\alpha(G)}$ has relative norm averaging (RNA) if for any finite subset $\mathcal{F}\subseteq D\rtimes_r G$ and $\varepsilon>0$, there exists a unitary $w\in (D\rtimes_r G) \cap \left(D^G\right)'$ such that
    \begin{align}\label{relative norm averaging}
       \| E(w^*xw)-P_K(E(x)) \|\leq \varepsilon
    \end{align}
    for every $x\in \mathcal{F}$.
\end{dfn}

Thus the Haar average of the identity Fourier coefficient is approximated, simultaneously on finite sets, by conjugation with a single unitary followed by $E$. Requiring $w\in (D^G)'$ ensures that this conjugation fixes $D^G$ pointwise.

\begin{rmk}\label{norm averaging remark}
We use the following norm averaging property (NA), introduced
in \cite{BellSelflessNonnuclear} in the tracial setting.
A $\sigma$-preserving action
$\alpha:G\to\aut(D,\sigma)$ is said to have NA if for any finite subset $\mathcal{F}\subseteq D\rtimes_r G$ and $\varepsilon>0$, there exists a unitary $w\in D\rtimes_r G$ such that
    \begin{align}\label{norm averaging}
       \| E(w^*xw)-\sigma(E(x)) \|\leq \varepsilon
    \end{align}
    for every $x\in \mathcal{F}$.

    Suppose that $\alpha$ is almost periodic and ergodic, in the sense that
$D^G=\mathbb C1_D$. Then since
    \begin{align*}
        P_K(a)=\sigma(P_K(a))=\int_{K}\sigma(\kappa(a)) dm_K(\kappa)=\int_{K}\sigma(a) dm_K(\kappa)=\sigma(a),
    \end{align*}
and the commutation requirement in RNA is automatic, RNA coincides with NA in this case. Moreover, $P_K$ is the
unique $K$-invariant state, which is tracial by \cite{StormerCompactErgodic}. 
\end{rmk} 

\begin{rmk}\label{tracial state on crossed}
    If $\sigma$ is a $G$-invariant tracial state on $D$, then $\sigma\circ E$ is again tracial on $D\rtimes_r G$. Therefore, if $\alpha$ is the almost periodic ergodic action, then $P_K\circ E$ is a tracial state.
\end{rmk}

Our first main result establishes RNA under a connectedness assumption.

\begin{thm}\label{main 1}
    Let $D$ be a unital separable $C^*$-algebra with a faithful state $\sigma$ and $G$ be a countable discrete group with an almost periodic $\sigma$-preserving action $\alpha:G \to \aut(D,\sigma)$.
    We assume that $K:=\overline{\alpha(G)}\subseteq \aut(D)$ is connected and nontrivial. 
    Then the action $\alpha$ has RNA.
\end{thm}

We recall the notions of selflessness and complete selflessness;
see \cite{RobertSelfless} and \cite[Section~6]{OzawaProximality}.

\begin{dfn}
Let $A$ be a separable unital $C^*$-algebra with GNS-faithful state $\varphi$.
We call $(A,\varphi)$ selfless if there exist a unital
$C^*$-algebra $Q\neq\mathbb C$ with a GNS-faithful state $\psi$,
a free ultrafilter $\omega$, and a unital state preserving embedding
\[
    \Theta:(A,\varphi)*(Q,\psi)\hookrightarrow A^\omega,
\]
such that $\Theta|_A=\Delta_A$, where $\Delta_A:A\hookrightarrow A^{\omega}$ is the diagonal embedding.

The pair $(A,\varphi)$ is \emph{completely selfless} if
$\Theta$ can be chosen so that, for every $C^*$-algebra $C$,
the induced homomorphism $((A,\varphi)*(Q,\psi))\otimes C\to (A\otimes C)^{\omega}$ is continuous with respect to the minimal tensor norm.
\end{dfn}

To obtain selflessness, we combine this result with Bell's criterion below. Although \cite[Theorem~5.1]{BellSelflessNonnuclear} is stated for tracial states, its proof does not use traciality and applies to faithful invariant states.

\begin{thm}[\cite{BellSelflessNonnuclear} Theorem~5.1]\label{Bell main}
    Let $D$ be a unital separable $C^*$-algebra with a faithful state $\sigma$, $G$ be a countable discrete group with property PHP. If $\alpha:G\to \aut(D,\sigma)$ is an action with NA, then $(D\rtimes_r G, \sigma\circ E)$ is completely selfless.
\end{thm}

Therefore, combining Theorem~\ref{main 1} and Theorem~\ref{Bell main}, we obtain the following corollary.

\begin{cor}\label{cor selfless}
    Let $D$ be a unital separable $C^*$-algebra with a faithful state $\sigma$, $G$ be a countable discrete group with property PHP, which acts on $(D,\sigma)$ ergodically by a $\sigma$-preserving action $\alpha$. If $\overline{\alpha(G)}$ is a connected and compact subgroup of $\aut(D,\sigma)$, then the reduced crossed product $D\rtimes_r G$ is completely selfless, especially has stable rank one, strict comparison and unique trace. 
\end{cor}

Selflessness for almost periodic action was also concerned in \cite{OhshimaSelfless}, where the author assumed that a coefficient algebra $D$ is separable, simple, exact, $\mathcal{Z}$-stable, and monotracial. In our setting, we only assume that $D$ is unital and separable. In particular,
our results apply to commutative coefficient algebras arising from
compact homogeneous spaces.

\begin{ex}\label{example}
    Isometric actions on compact Riemannian manifolds give almost periodic
actions on their algebras of continuous functions. The invariant state
in the following examples is integration against normalized Riemannian
volume.
    \begin{enumerate}
        \item $G=\FF_2$ admits a dense embedding into $SO(3)$ \cite{SwierczkowskiRotation}, so we have a natural ergodic action $\alpha:G\to \aut(C(S^{2}))$ with $K=\overline{\alpha(G)}\cong SO(3)$. Since $\FF_2$ has PHP, the reduced crossed product $C(S^2)\rtimes_r G$ is completely selfless.
    
        \item Let $n\geq 3$ be an odd number. $G=SO\left(n;\ZZ[\frac{1}{5}]\right)$ is a dense subgroup of $SO(n;\RR)$, so we have a natural ergodic action $\alpha:G\to \aut(C(S^{n-1}))$ with $K=\overline{\alpha(G)}\cong SO(n;\RR)$.
        By \cite[Lemma~5.2]{VigdorovichSelfless}, $G$ has trivial amenable radical.
        Since a linear group with trivial amenable radical has PHP \cite[Proof of Theorem~1.4]{VigdorovichSelfless}, $SO\left(n;\ZZ[\frac{1}{5}]\right)$ has PHP if $n$ is odd. Therefore the reduced crossed product $C(S^{n-1})\rtimes_r G$ is completely selfless for odd $n\geq 3$. (See Example~\ref{example 2} for even $n$.)

        Note that since $SO\left(n;\ZZ[\frac{1}{5}]\right)$ has Kazhdan's property (T) for $n\geq 5$ and the action is ergodic \cite{Margulis1980Some}, the action $\alpha$ has a spectral gap. Therefore, the action is not squarely divisible \cite[Proposition~3.15]{BellGeffenKerrStable}
    \end{enumerate}
\end{ex}

\vspace{.2in}

We next study operator-valued selflessness when the fixed-point
algebra $D^G$ is nontrivial.
Let $A$ be a unital $C^*$-algebra and $B$ its $C^*$-subalgebra with an expectation $\Phi:A\to B$. We assume that $\Phi$ is GNS-faithful, meaning that the associated representation of $A$ on the Hilbert $B$-module
$L^2(A,\Phi)$ is faithful.

\begin{dfn}[\cite{GaoJungeOperatorValued} Definition~7.1]
    We say that the triple $(A, \Phi,B)$ is selfless if there exists a $C^*$-probability space $(C,\rho)$ with $C\neq\mathbb C$ and
GNS-faithful $\rho$, a free ultrafilter $\omega$ and unital embedding
    \begin{align*}
        \Theta:(A,\Phi)*_B(B\otimes C, \id_B\otimes \rho) \hookrightarrow (A,\Phi)^{\omega}
    \end{align*}
    satisfying that $\Theta|_A=\Delta|_A$ and $\Phi^{\omega}\circ \Theta=\Delta_A \circ (\Phi*_B (\id_B\otimes \rho))$.
\end{dfn}

\begin{rmk}
    The expectation $P_K:D\to D^G$ is faithful and so is $P_K\circ E:D\rtimes_r G \to D^G$.
\end{rmk}

The following is the main theorem of the second half of the paper, which can be thought as the relative generalization of \cite[Theorem~5.1]{BellSelflessNonnuclear}.

\begin{thm}\label{main 2}
Let $D$ be a unital separable $C^*$-algebra with a faithful state $\sigma$, $G$ be a countable discrete group which acts on $D$ by an almost periodic action.
    If $\alpha:G \to \aut(D,\sigma)$ has RNA and $G$ has PHP, then the triple $(D\rtimes_r G, P_K\circ E, D^G)$ is selfless.
\end{thm}

\section{Relative norm averaging property for almost periodic actions}
In this section, we prove Theorem~\ref{main 1}. We provide analytic lemmas and a group theoretic lemma in separated subsections and then complete the proof in the third subsection. We continue assuming that $D$ is a unital separable $C^*$-algebra and $\sigma$ is a faithful state on $D$.

\subsection{Analytic ingredient}
For a fixed parameter $t\in (0,1]$, we define a function $\varphi(z):=\exp(it(z+\overline{z}))$ for $z\in \mathbb{T}\subseteq \CC$. Let $\varphi(z)=\sum_{n\in \ZZ}c_n z^n$ be the Laurent series of $\varphi$. Then we have the following
\begin{enumerate}[i)]
    \item $\displaystyle
    \sum_{n\in \ZZ} |c_n|
    \leq \sum_n
    \sum_{\substack{a-b=n\\ a,b\geq 0}}
    \frac{t^{a+b}}{a!b!}
    \leq \sum_{a,b\geq 0}
    \frac{t^{a+b}}{a!b!}
    = \exp(2t)<\infty$

    \item $\sum_{n\in \ZZ} |c_n|^2=\|\varphi\|^2_2=1$

    \item $c_0\neq 0$, $c_1\neq 0$.
\end{enumerate}
We denote by $p_n:=|c_n|^2$ and for each $\kappa\in \aut(D,\sigma)$, we define a $\kappa$-average by
\begin{align*}
    R_{\kappa}(d)=\sum p_n\kappa^{-n}(d).
\end{align*}
It is easy to see the following.
\begin{lem}
Let $K\subseteq \aut(D,\sigma)$ be a compact subgroup.
    For any $\kappa\in K$, $R_{\kappa}$ is contractive with respect to both $C^*$-norm and $L^2$-norm and satisfies $$P_KR_{\kappa}=R_{\kappa}P_K=P_K.$$
\end{lem}

\begin{lem}\label{l2-average}
    For any $d\in D$ with $P_K(d)=0$, we have that
    \begin{align*}
        \int_K\|R_{\kappa}(d)\|_{L^2(D,\sigma)}^2dm_K(\kappa)\leq (1-2p_0p_1)\|d\|_{L^2(D,\sigma)}^2.
    \end{align*}
\end{lem}

\begin{proof}
    First, note that
    \begin{align*}
        \|R_{\kappa}(d)\|_{L^2(D,\sigma)}^2=&\sum_{n,m\in \ZZ}p_np_m\left\langle \kappa^{-n}(d),\kappa^{-m}(d)\right\rangle_{L^2(D,\sigma)}\\
        =&\frac{1}{2}\left(\sum_{n,m\in \ZZ}p_np_m\left(\|\kappa^{-n}(d)\|_{L^2(D,\sigma)}^2+\|\kappa^{-m}(d)\|_{L^2(D,\sigma)}^2 \right)\right.\\
        &-\left.\sum_{n,m\in \ZZ} p_np_m\|\kappa^{-n}(d)-\kappa^{-m}(d)\|_{L^2(D,\sigma)}^2\right)\\
        =&\|d\|^2_{L^2(D,\sigma)}-\frac{1}{2}\sum_{n,m\in \ZZ} p_np_m\|\kappa^{-n}(d)-\kappa^{-m}(d)\|_{L^2(D,\sigma)}^2.
    \end{align*}
    Therefore,
    \begin{align}\label{before int}
        \|d\|^2_{L^2(D,\sigma)}-\|R_{\kappa}(d)\|_{L^2(D,\sigma)}^2&=\frac{1}{2}\sum_{n,m\in \ZZ} p_np_m\|\kappa^{-n}(d)-\kappa^{-m}(d)\|_{L^2(D,\sigma)}^2\\
        &\geq p_0p_1\|d-\kappa^{-1}(d)\|_{L^2(D,\sigma)}^2
    \end{align}
    But since $P_K(d)=0$, we have
    \begin{align*}
        \int_K \|d-\kappa^{-1}(d)\|_{L^2(D,\sigma)}^2 dm_K(\kappa)=&2\int_K \|d\|_{L^2(D,\sigma)}^2 dm_K(\kappa)\\
        &-2 \Re\left(\left\langle d,\int_K\kappa^{-1}(d)dm_K(\kappa)\right\rangle_{L^2(D,\sigma)}\right)\\
        =&2\int_K \|d\|_{L^2(D,\sigma)}^2 dm_K(\kappa).
    \end{align*}
    Therefore, by integrating \eqref{before int}, we obtain the desired inequality.
\end{proof}

\begin{lem}\label{average}
    If $K\subseteq \aut(D,\sigma)$ is a compact subgroup and $d_1,\cdots, d_r\in D$ satisfy $P_K(d_j)=0$ for all $j=1,\cdots,r$, then for every $\varepsilon>0$, there exists $\kappa_1, \cdots, \kappa_k\in K$ such that
    \begin{align*}
        \|R_{\kappa_k}\circ \cdots \circ R_{\kappa_1}(d_j)\|\leq \varepsilon
    \end{align*}
    for all $j=1,\cdots,r$.
\end{lem}

\begin{proof}
We may assume that $\|d_j\|_{L^2(D,\sigma)}=1$ for all $j=1,\cdots,r$.

Since $K\subseteq \aut(D,\sigma)$ is a compact subset, for each $d\in D$, $\{\kappa(d):\kappa\in K\}$ is compact. Therefore the set
\begin{align*}
    C:=\bigcup_{j=1}^r \overline{\conv}(\{\kappa(d_j):\kappa\in K\})
\end{align*}
is again compact. Since $\sigma$ is faithful, for a given $\varepsilon>0$, there exists $\delta>0$ such that
\begin{align*}
    \|d\|_{L^2(D,\sigma)}\leq \delta \Longrightarrow \|d\|\leq \varepsilon
\end{align*}
for every $d\in C$.
    Let $d_1,\cdots, d_r\in D$ with $P_K(d_j)=0$ be given. By Lemma~\ref{l2-average}, we have
    \begin{align*}
        \int_K\sum_{j=1}^r\|R_{\kappa}(d_j)\|_{L^2(D,\sigma)}^2dm_K(\kappa)\leq (1-2p_0p_1)\sum_{j=1}^r\|d_j\|_{L^2(D,\sigma)}^2.
    \end{align*}
    Therefore, for some $\kappa_1$, we have that
    \begin{align*}
        \sum_{j=1}^r\|R_{\kappa_1}(d_j)\|_{L^2(D,\sigma)}^2\leq (1-2p_0p_1)\sum_{j=1}^r\|d_j\|_{L^2(D,\sigma)}^2.
    \end{align*}
    By repeating this, one can find a positive integer $k$ such that for every $j=1,\cdots ,r$,
    \begin{align*}
        \|R_{\kappa_k}\circ \cdots \circ R_{\kappa_1}(d_j)\|_{L^2(D,\sigma)}^2&\leq \sum_{j=1}^r\|R_{\kappa_k}\circ \cdots \circ R_{\kappa_1}(d_j)\|_{L^2(D,\sigma)}^2\\
        &\leq \frac{\delta^2}{r}\sum_{j=1}^r\|d_j\|_{L^2(D,\sigma)}^2
        \leq \delta^2.
    \end{align*}
    Therefore $\|R_{\kappa_k}\circ \cdots \circ R_{\kappa_1}(d_j)\|\leq \varepsilon$.
\end{proof}

\subsection{Group theoretic ingredient}
We need the following lemma to avoid nontrivial cancellations in the computation \eqref{computation of expectation}.
\begin{lem}\label{non-cancellation}
    Let $\alpha:G \to \aut(D,\sigma)$ be a nontrivial almost periodic action with $K=\overline{\alpha(G)}$ connected. Then for any positive integer $L$, family of nonempty open subsets $V_1, \cdots, V_k$ in $K$ and a finite subset $F\subseteq G$, there exist $h_1, \cdots h_k\in G$ with $\alpha_{h_j}\in V_j$ such that for any $(n_1,\cdots,n_k), (m_1,\cdots,m_k)\in \{-L,\cdots,L\}^k$ and $g\in F\setminus \{e\}$, we have
    \begin{align}
        h_1^{n_1}\cdots h_k^{n_k}=h_1^{m_1}\cdots h_k^{m_k} \Longrightarrow n_j=m_j \quad\forall j=1,\cdots ,k\\
        (h_1^{n_1}\cdots h_k^{n_k})^{-1}g(h_1^{m_1}\cdots h_k^{m_k})\neq e.
    \end{align}
\end{lem}

\begin{proof}
\noindent{\textbf{Case 1 $K$:Lie group}}

First, we assume that $K$ is a compact connected non-trivial Lie group. For distinct $\vec{n},\vec{m}\in \{-L,\cdots,L\}^k$, we define a real analytic non-constant function
\begin{align}\label{analytic function}
    W_{\vec{n},\vec{m}}:K^k\to K,\; (\kappa_1,\cdots, \kappa_k)\mapsto (\kappa_1^{n_1}\cdots\kappa_k^{n_k})(\kappa_1^{m_1}\cdots\kappa_k^{m_k})^{-1}.
\end{align}
Since $K$ is connected, $W_{\vec{n},\vec{m}}^{-1}(\kappa)$ has no interior for any $\kappa\in K$, because if not, $W_{\vec{n},\vec{m}}$ would be constant by the connectedness of $K$ and $W_{\vec{n},\vec{m}}$ being real analytic. Therefore
\begin{align*}
M:=\left\{
(\kappa_1,\ldots,\kappa_k)\in K^k:
\begin{aligned}
& W_{\vec n,\vec m}((\kappa_1,\ldots,\kappa_k))
  \notin \{\alpha_g:g\in F\cup\{e\}\}\\
& \text{for all }\vec n,\vec m\in\{-L,\ldots,L\}^k
  \text{ with }\vec n\neq\vec m
\end{aligned}
\right\}.
\end{align*}
is an open dense subset in $K^k$ and $M\cap(V_1\times \cdots \times V_k)$ is a nonempty open set. Since $\alpha(G)$ is dense in $K$, we can find $(h_1,\cdots,h_k)\in G^k$ such that 
\begin{align*}
    (\alpha(h_1),\cdots,\alpha(h_k))\in M\cap(V_1\times \cdots \times V_k).
\end{align*}
By the definition of $W_{\vec{n},\vec{m}}$, if $\vec{n}\neq \vec{m}$, then $\alpha_{h_1^{n_1}\cdots h_k^{n_k}}\neq \alpha_g\cdot\alpha_{h_1^{m_1}\cdots h_k^{m_k}}$ and so $h_1^{n_1}\cdots h_k^{n_k}\neq g (h_1^{m_1}\cdots h_k^{m_k})$ for all $g\in F$. If $\vec{n}= \vec{m}$, then we automatically have $h_1^{n_1}\cdots h_k^{n_k}\neq g (h_1^{m_1}\cdots h_k^{m_k})$ for all $g\in F\setminus \{e\}$.

\noindent{\textbf{Case 2 general case}}

Next, we assume that $K$ is a general compact connected group. By Peter-Weyl Theorem, there exists a non-trivial compact connected Lie group $K'$ and a surjective continuous homomorphism $\pi:K\twoheadrightarrow K'$. Note that since $\pi$ is an open map, $\pi(V_j)$ is a nonempty open set in $K'$. By applying the same argument as the previous Lie group case to find an open dense subset $M'$ of $(K')^k$ such that for any points $(\kappa_1,\cdots, \kappa_k)\in M'$ and $g\in F\setminus \{e\}$, we have that
\begin{align*}
    \kappa_1^{n_1}\cdots \kappa_k^{n_k}=\kappa_1^{m_1}\cdots \kappa_k^{m_k} \Longrightarrow n_j=m_j \quad\forall j=1,\cdots ,k\\
        (\kappa_1^{n_1}\cdots \kappa_k^{n_k})^{-1}\pi(\alpha(g))(\kappa_1^{m_1}\cdots \kappa_k^{m_k})\neq e \quad  \forall \vec{n}\neq\vec{m}\in \{-L,\cdots,L\}^k.
\end{align*}
Since $(\pi^k)^{-1}(M')$ is dense in $K^k$, it has a nontrivial intersection with $V_1\times \cdots \times V_k$. By taking $(h_1,\cdots,h_k)\in G^k$ in a way that 
\begin{align*}
    (\alpha(h_1),\cdots,\alpha(h_k))\in (\pi^k){-1}(M')\cap(V_1\times \cdots \times V_k)
\end{align*}
we can complete the proof similarly as the previous case.
\end{proof}

\subsection{Proof of Theorem~\ref{main 1}}

In this subsection, we complete the proof of Theorem~\ref{main 1}.

Let $\{x_l\}_{l=1}^r=\mathcal{F}\subseteq D\rtimes_r G$ and $\varepsilon>0$ be given. We denote by $y_l:=x_l-P_K(E(x_l))$. Then since $P_K(E(x_l))\in D^G$, we have 
$E(w^*P_K(E(x_l))w)=P_K(E(x_l))$ for any unitary $w$ in $(D^G)'$. Therefore, it suffices to show that there exists a unitary $w\in C^*_r(G)$ such that 
\begin{align}\label{y-estimate}
    \|E(w^*y_lw)\|\leq \varepsilon.
\end{align}

We may assume that $y_l\in C_c(G,D)$ and denote $y_l=\sum_{g \in F} d_{l,g}g$ for a common finite subset $F\subseteq G$ and set $b_l:=d_{l,e}$. Note that $P_K(E(b_l))=0$. By Lemma~\ref{average}, there exist a positive integer $k$ and $\kappa_1',\cdots, \kappa_k'\in \overline{\alpha(G)}$ such that
\begin{align*}
    \|R_{\kappa_k'}\circ \cdots \circ R_{\kappa_1'}(b_l)\| \leq \frac{\varepsilon}{6}.
\end{align*}
For each $j$, we can take an open neighborhood $V_j\subseteq \aut(D,\sigma)$ such that for any $\kappa_1\in V_1, \cdots, \kappa_k\in V_k$,
\begin{align*}
    \|R_{\kappa_k}\circ \cdots \circ R_{\kappa_1}(b_l)\| \leq \frac{\varepsilon}{3}.
\end{align*}

For positive integer $L$, we define $q_L:=\sum_{|n|\leq L} c_nz^n$. Fix enough large $L$ such that
\begin{align*}
     k 2^k\|\varphi-q_L\|_{C(\mathbb T)} &\leq 
    \frac{\varepsilon}{9\max\{\|y_1\|, \cdots, \|y_r\|\}}\\
    \|q_L\|^k&\leq 2\\
  k 2^k \cdot \sum_{|n|>L} p_n &\leq \frac{\varepsilon}{6\max\{\|b_1\|, \cdots, \|b_r\|\}}.  
\end{align*}
For this $L$, finite subset $F\subseteq G$ and open sets $V_1,\cdots,V_k$, we take $h_j\in G$ as in Lemma~\ref{non-cancellation}
and define $w=\prod_{j=1}^k \varphi(h_j)$. We will show that $w$ satisfies the inequality \eqref{y-estimate}. Put $Q_L:=\prod_{j=1}^k q_L(h_j)$.
We denote for $\Vec{n}\in \{-L,\cdots,L\}^k$, $c_{\Vec{n}}:=c_{n_1}\cdots c_{n_k}\in \CC$ and $h_{\Vec{n}}:=h_1^{n_1}\cdots h_{k}^{n_k}\in G$. Then by the non-cancellation criteria given in Lemma~\ref{non-cancellation},
\begin{align}\label{computation of expectation}
    E(Q_L^*y_lQ_L)&=E\left(\sum_{\Vec{n},\Vec{m}\in \{-L,\cdots,L\}^k} \overline{c_{\Vec{n}}}c_{\Vec{m}} h_{\Vec{n}}^*\left(\sum_{g\in F}d_{l,g}g\right) h_{\Vec{m}}\right)\\
    &=E\left(\sum_{\Vec{n},\Vec{m}\in \{-L,\cdots,L\}^k} \sum_{g\in F}\left(\overline{c_{\Vec{n}}}c_{\Vec{m}} \alpha_{h_{\Vec{n}}}^{-1}(d_{l,g})\right) h_{\Vec{n}}^*gh_{\Vec{m}} \right)\\
    &=\sum_{\Vec{n}\in \{-L,\cdots,L\}^k} |c_{\Vec{n}}|^2 \alpha_{h_{\Vec{n}}}^{-1}(b_l)\\
    &=\left( \sum_{n_k=-L}^L p_{n_k} \alpha_{h_k}^{-n_k} \cdots \left( \sum_{n_2=-L}^L p_{n_2} \alpha_{h_2}^{-n_2} \left(\sum_{n_1=-L}^Lp_{n_1}\alpha_{h_1}^{-n_1}(b_l)\right) \right)\cdots \right).
\end{align}
Now if we denote by $R_{\kappa,L}:=\sum_{n=-L}^L p_n\kappa^{-n}:D\to D$, then 
\begin{align*}
    E(Q_L^*y_lQ_L)=R_{\alpha_{h_k},L}\circ \cdots \circ R_{\alpha_{h_1},L}(b_l).
\end{align*}
Also for any $\kappa\in K$, we have $\|R_{\kappa,L}(b_l)-R_{\kappa}(b_l)\|\leq \frac{\varepsilon}{6k2^k\max\{\|b_1\|, \cdots, \|b_r\|\}}\|b_l\|$ and so 
\begin{align*}
    \|E(Q_L^*y_lQ_L)\|\leq\| R_{\alpha_{h_k}}\circ \cdots \circ R_{\alpha_{h_1}}(b_l)\|+\frac{\varepsilon}{3}\leq \frac{2\varepsilon}{3}.
\end{align*}
Since $ \|w-Q_L\|\leq k2^k\|\varphi-q_L\|\leq\frac{\varepsilon}{9}$, we have
\begin{align*}
    \|E(w^*y_lw)\|&\leq \|E(w^*y_lw)-E(Q_L^*y_lQ_L)\|+\|E(Q_L^*y_lQ_L)\|\\
    &\leq \|y_l\|\cdot(1+\|Q_L\|)\cdot\|w-Q_L\|+\|E(Q_L^*y_lQ_L)\|\\
    &\leq \frac{(1+2)\varepsilon}{9}+\frac{2\varepsilon}{3}=\varepsilon.
\end{align*}
\qed

\section{Operator-valued selflessness for non-ergodic actions}

In this section, we prove Theorem~\ref{main 2}. First, we recall the definition of property PHP.

\begin{dfn}[\cite{OzawaProximality}]
    We say a discrete group $G$ has property PHP if for every $\varepsilon>0$ and finite subset $F\subseteq G$, there are a positive integer $k$, subsets $C_i\subseteq D_i \subseteq G$ and $t_i\in G$ for $i=1,\cdots ,k$ such that members 
    \begin{align*}
        \{aC_i:a\in F, i=1,\cdots,k\}\cup \{bt^{-1}_j(G\setminus D_i):b\in F, j=1,\cdots,k\}
    \end{align*}
    are pairwise disjoint and we have
    \begin{align*}
        \sup_{g\in G}|\{i:g\in D_i\}\cup \{j:g\in t^{-1}_j(G\setminus C_i)\}| \leq \varepsilon\sqrt{k}
    \end{align*}
\end{dfn}

\begin{proof}[Proof of Theorem~\ref{main 2}]

We consider the regular representation of $D\rtimes_r G$ on the Hilbert $D$-module $\ell^2(G)\otimes D$. Let $B:=\mathcal{L}_D(\ell^2(G)\otimes D)$ be the algebra of all adjointable operators on $\ell^2(G)\otimes D$. Note that $D\rtimes_r G$ is the norm completion of $C_c(G,D)$ in $B$. For any subset $I\subseteq G$, we denote by $\chi_I\in B(\ell^2(G))$ the projection onto $\ell^2(I)$.
The canonical expectation $E:D\rtimes_r G\to D$ can be extended to $\widetilde{E}:T\mapsto \langle (\delta_{e}\otimes 1_D),T(\delta_{e}\otimes 1_D)\rangle$ on $B$. Therefore $P_K\circ \widetilde{E}:B\to D^G$ is a ucp map extending $P_K\circ E$.

By \cite[Proposition~5.9,Theorem~7.3]{GaoJungeOperatorValued}, it suffices to show that for every finite subset $F$ of the unit ball of $D\rtimes_r G$ and $\varepsilon>0$, there exists $v\in B$ such that
\begin{align*}
    \|v^*xv-P_K\circ E(x)\|\leq \varepsilon\\
    \|P_K\circ \widetilde{E}(xvv^*x^*)\|\leq \varepsilon\\
    d(v+v^*, D\rtimes_r G)\leq \varepsilon
\end{align*}
for every $x\in F$.

Take $w\in (D\rtimes _r G)\cap (D^G)'$ in the definition of RNA for $F$ and $\frac{\varepsilon}{4}$.
For each $x\in F$, we take $y_x\in C_c(G,D)$ satisfying $\|w^*xw-y_x\|\leq \frac{\varepsilon}{4}$ and take a finite subset $S\subseteq G$, which contains the support of all $y_x$. Replacing $S$ by a larger finite subset, we may assume that
\begin{align*}
    \|(1-\chi_S) w^*x^*(\delta_e\otimes 1_D)\|\leq \frac{\varepsilon}{2}.
\end{align*}
for every $x\in F$.

Since $G$ has PHP, by the proof of \cite[Theorem~14]{OzawaProximality}, there are $n\geq \frac{4|S|}{\varepsilon^2}$ and a contraction $T\in B(\ell^2(G))$ of the form
\begin{align*}
    T:=\frac{1}{\sqrt{2n}}\sum_{i=1}^n (\chi_{C_i}t_i+t_i^{-1}(1-\chi_{D_i})),
\end{align*}
which satisfies
\begin{enumerate}[i)]
    \item $\|T^*T-1\|\leq \frac{\varepsilon}{4}$,

    \item $d(T+T^*,C^*_r(G))\leq \varepsilon$,

    \item $T^*(dg)T=0$ for every $d\in D$ and $g\in S\setminus \{e\}$,

    \item $\|\chi_S T\|^2\leq \frac{|S|}{n}$ for every finite subset $S\subseteq G$.
\end{enumerate}

 Define $v=w(T\otimes 1_D)w^*$.
 Since $w$ is a unitary and belongs to $D\rtimes _r G$, i) and ii) imply that $\|v^*v-1\|\leq \varepsilon$ and $d(v+v^*, D\rtimes_r G)\leq \varepsilon$. 
 By iii), we have that
\begin{align*}
    (T^*\otimes 1_D)y_x(T\otimes 1_D)=(T^*\otimes 1_D)E(y_x)(T\otimes 1_D).
\end{align*} 
Therefore, using $w, T\otimes 1\in (D^G)'$, we obtain
\begin{align*}
    &\|v^*xv-P_K(E(x))\|\\
    \leq &\|v^*xv-w(T^*\otimes 1)y_x(T\otimes 1)w^*\|+\|w(T^*\otimes 1)E(y_x)(T\otimes 1)w^*-P_K(E(x))\|\\
    \leq &\|w^*xw-y_x\|+\|(T^*\otimes 1)E(y_x)(T\otimes 1)-P_K(E(x))\|\\
    \leq & \|w^*xw-y_x\|+\|(T^*\otimes 1)E(y_x)(T\otimes 1)-(T^*\otimes 1)E(w^*xw)(T\otimes 1)\|\\
    &+\|(T^*\otimes 1)E(w^*xw)(T\otimes 1)-(T^*\otimes 1)P_K(E(x))(T\otimes 1)\|\\
    &+\|(T^*\otimes 1)(T\otimes 1)P_K(E(x))-P_K(E(x))\|\\
    \leq & \varepsilon.
\end{align*}
Also, by iv) and the choice of $n\geq \frac{4|S|}{\varepsilon^2}$, we have
\begin{align*}
     \|P_K\circ \widetilde{E}(xvv^*x^*)\|^{\frac{1}{2}}&\leq \|\widetilde{E}(xvv^*x^*)\|^{\frac{1}{2}}\\
     &= \|(T^*\otimes 1)w^*x^*(\delta_{e}\otimes 1_D)\| \\
     &\leq \|(T^*\otimes 1)\chi_S w^*x^*(\delta_e\otimes 1_D)\|+\|(T^*\otimes 1)(1-\chi_S) w^*x^*(\delta_e\otimes 1_D)\|\\
     &\leq \sqrt{\frac{|S|}{n}}\|x\|+\frac{\varepsilon}{2}\leq \varepsilon.
\end{align*}

\end{proof}

\section{Selflessness of crossed products by groups with the finite radical}

In this section, we extend Corollary~\ref{cor selfless} to the case where $G$ has finite radical, which covers the second case of Example~\ref{example} for even $n\geq 3$.

\begin{cor}\label{finite radical}
    Let $D$ be a unital separable $C^*$-algebra with a faithful state $\sigma$, $G$ be a countable discrete group that acts on $(D,\sigma)$ by a $\sigma$-preserving action $\alpha$ and $N$ be a finite normal subgroup of $G$. We assume the following
    \begin{enumerate}
        \item the quotient group $Q:=G/N$ satisfies PHP,
        \item  the action $\beta:Q\curvearrowright D^N$ induced by $\alpha$ is ergodic and almost periodic and $\overline{\beta(Q)}$ is connected and nontrivial,
        \item the projection $p:=\displaystyle \frac{1}{|N|}\sum_{g\in N} g$ is full in $D\rtimes_r G$.
    \end{enumerate}
    In this case, $(D\rtimes_r G, \sigma\circ E)$ is completely selfless.
\end{cor}
For tracial $C^*$-probability spaces, complete selflessness of the full corner is equivalent to that of the original $C^*$-algebra \cite[Lemma~5.2]{KlissePagliero}. By Remark~\ref{norm averaging remark} and Remark~\ref{tracial state on crossed}, $\sigma\circ E$ is a tracial state on $D^N\rtimes_r Q$. Once we prove that $\sigma\circ E$ is also a trace on $D\rtimes_r G$, then the Corollary follows from Lemma~\ref{trace isom} using the aforementioned permanence property. The trace $\sigma\circ E$ on $p(D\rtimes_r G)p\cong D^N\rtimes_r Q$ extends to a trace $\tau$ of $D\rtimes_r G$ by \cite[Remark~3.11]{RordamFixedPoints}, which is continuous by the fullness of $p$. We will prove that $\tau=\sigma\circ E$ on $D$. On $D^N$, they agree since $\beta$ is ergodic and almost periodic and so $P_{\overline{\beta(Q)}}$ is the unique $Q$-invariant state (Remark~\ref{norm averaging remark}). Then by $G$-invariance we have that for all $d\in D$
\begin{align*}
    \tau|_D(d)=\tau|_D\left(\frac{1}{|N|} \sum_{g\in N}\alpha_g(d)\right)=\sigma\left(\frac{1}{|N|} \sum_{g\in N}\alpha_g(d)\right)=\sigma(d).
\end{align*}
So by Remark~\ref{tracial state on crossed}, $\sigma \circ E$ is tracial on $D\rtimes_r G$.

\begin{lem}\label{trace isom}
    Under the assumption of Corollary~\ref{finite radical}, we have a trace preserving isomorphism;
    \begin{align*}
        \Phi:(D^N\rtimes_r Q,\sigma\circ E) \to \left(p(D\rtimes_r G)p, \frac{1}{\sigma\circ E(p)}\sigma\circ E\right)=\left(p(D\rtimes_r G)p, |N|\sigma\circ E\right).
    \end{align*}
    
\end{lem}

\begin{proof}
    First, we probe the first isomorphism. Since $N$ is normal, for $g\in G$ we have $pg=gp=pgp$ and for any $a\in D^N$ we have $ap=pa=pap$. Therefore 
    \begin{align*}
    \Phi:D^N\rtimes Q \to p(D\rtimes_r G)p;\; 
    \left\{
    \begin{array}{cc}
       a\mapsto pap  & \text{ (for $a\in D^N$)} \\
        {[g]} \mapsto pgp & \text{ (for $[g]\in Q$)}.
    \end{array}
    \right.
    \end{align*}
    is an well-defined surjective $*$-homomorphism on the algebraic crossed product. For a finite sum $a:=\sum_{[g]} a_{[g]} g\in D^N\rtimes_r Q$, $\sigma\circ E (a)=a_{[e]}$. On the other hand, we have
    \begin{align*}
        |N|\sigma\circ E (\Phi(a))=\sigma\circ E\left(\sum_{[g]}\sum_{n\in N} a_{[g]}gn\right)=\sigma (a_{[e]}).
    \end{align*}
    Therefore $\Phi$ is state preserving and which induces a unitary
    \begin{align*}
        L^2(D^N\rtimes Q, \sigma\circ E) \to L^2(\left(p(D\rtimes_r G)p, |N|\sigma\circ E\right)).
    \end{align*}
    Since $D^N\rtimes Q$ and $p(D\rtimes_r G)p$ are faithfully represented on these Hilbert spaces, $\Phi$ preserves their $C^*$-norms.
\end{proof}

\begin{ex}\label{example 2}
    We extend the example of an action $ \alpha:SO\left(n;\ZZ[\frac{1}{5}]\right) \curvearrowright C(S^{n-1})$ for even $n\geq 3$. (See Example~\ref{example} for odd $n\geq 3$.) Let $G:=SO\left(n;\ZZ[\frac{1}{5}]\right)$ and $N:=\{-1,1\}$. Then since the quotient $Q=G/N$ can be faithfully represented on $\GL(M_n(\RR))$ by the adjoint representation with image dense in $PSO(n;\RR)$, it has trivial amenable radical \cite[Lemma~5.2]{VigdorovichSelfless} and thus has PHP \cite[Proof of Theorem~1.4]{VigdorovichSelfless}. The induced action $\beta:Q\curvearrowright C(\RR P^{n-1})$ is ergodic with $\overline{\beta(Q)}=PSO(n;\RR)$. Therefore it remains to prove that $$p=\frac{I_n+(-I_n)}{2}\in C(S^{n-1})\rtimes_r G$$ is full. Since $N$ acts freely on $S^{n-1}$, there exists a finite cover $\{U_1,\cdots,U_l\}$ of $S^{n-1}$ such that $U_i\cap (-U_i)=\emptyset$ for all $i=1,\cdots,l$. Let $\{\mu_i\}_{i=1}^l$ be a partition of unity subordinate to $\{U_i\}_{i=1}^l$. Then since $\mu_i\alpha_{-I_n}(\mu_i)=0$, we have that
    \begin{align*}
        \sum_i \sqrt{\mu_i}p\sqrt{\mu_i}=\frac{1}{2}\sum_i \sqrt{\mu_i}(\sqrt{\mu_i}+\alpha_{-I_n}(\sqrt{\mu_i}))=\frac{I_n}{2}.
    \end{align*}
    Therefore, $p$ is full. By Corollary~\ref{finite radical}, $ C(S^{n-1})\rtimes_r SO\left(n;\ZZ[\frac{1}{5}]\right)$ is completely selfless for all $n\geq 3$.
\end{ex}

\section*{Acknowledgements}
The author would like to thank Mikael R\o rdam for discussions on the topic.

\section*{AI Statement}
The proofs were developed with the help of Astra (OpenAI) and Opus 5.5 (Claude), which were also used to improve the writing and presentation. The author has checked all proofs and takes full responsibility for the final manuscript.

\bibliographystyle{alpha}
\bibliography{main}

\end{document}